\documentclass[12pt]{amsart}
\usepackage{amscd,amsthm,amsmath} 
\usepackage{pb-diagram}
\newtheorem{thm}{Theorem}[section]

\newtheorem{lem}[thm]{Lemma}
\newtheorem{prop}[thm]{Proposition}

\newtheorem{claim}[thm]{Claim}
\theoremstyle{definition}
\newtheorem{defn}[thm]{Definition}
\theoremstyle{remark}
\newtheorem{rem}[thm]{Remark}
\numberwithin{equation}{section}

\begin{document}

\title{Mislin genus of  maps}%
\author{Jianzhong Pan }%
\address{Institute of Math.,Academia Sinica ,Beijing China and
 Department of Mathematics Education , Korea University , Seoul , Korea }%
\email{pjz@math03.math.ac.cn}%
\author{Moo Ha Woo}
\address{ Department of Mathematics Education , Korea University , Seoul , Korea}
\email{woomh@kuccnx.korea.ac.kr}
\thanks{The first author is partially supported by the NSFC
projects 19701032 , 10071087 and ZD9603 of Chinese Academy of
Science and Brain Pool program of KOSEF and the second author
wishes to acknowledge the financial supports of the Korea Research
Foundation made in the
program year of (2000)}%

\begin{abstract}
 In this paper, we prove that the Mislin genus of a (co-)H-map
 between (co-)H-spaces under certain natural conditions is a
 finite abelian group which generalizes  results in Zabrodsky\cite{az1} , McGibbon\cite{mcg} and
 Hurvitz\cite{hurv}
\end{abstract}
\maketitle
 \section{Introduction} \label{S:intro}
  Let $X$ be a nilpotent connected CW complex of finite type and
  $X_{(p)}$ be the $p$-localization. The Mislin genus\cite{mislin} of
  $X$ is the set $G(X)$ of the homotopy types of nilpotent connected CW
  complexes Y of finite type such that $Y_{(p)}\simeq X_{(p)}$ for all
  prime p. In general $G(X)$ is not trivial. The first general
  result about Mislin genus is that of Wilkerson:
  \begin{thm}\cite{wilk}
  Let $X$ be a 1-connected CW-complex of finite type. If
  $H_n(X,\mathbb{Z})=0$ for $n$ sufficiently large  or $\pi_n(X)=0$ for $n$
  sufficiently large , then $G(X)$ is finite.
  \end{thm}
  In general it is difficult to determine the set $G(X)$. An
  exceptional case is the following result of
  Zabrodsky\cite{az1}\cite{az2}.
  \begin{thm} Let $X$ be a
  1-connected rational H-space (i.e., space of finite type such that
  its rationalization is an H-space)  with only finitely number
  of nontrivial homotopy groups. Then the following is an exact
  sequence
  $$[X,X]_{\hat{t}} \overset{d}{\rightarrow}(\mathbb{Z}^*_{\hat{t}}/\pm1)^l \overset{\xi}{\rightarrow}G(X) \rightarrow 0 $$
where $\hat{t}$ is a certain positive integer depending on $X$ ,
$l$ is the number of integers $k$ with $QH^k(X;\mathbb{Q})\neq0$
and $[X,X]_{\hat{t}}$ is the set of homotopy classes of the
self-maps of $X$ which are $\hat{t}$-equivalences.
  \end{thm}
\begin{rem}
 $QH^*(X;\mathbb{Q})$ in the above Theorem is the indecomposable
 module of the algebra  $H^*(X,\mathbb{Q})$ over $\mathbb{Q}$.
\end{rem}
Recently McGibbon \cite{mcg} generalized this to the case of
connected nilpotent rational H-space and to the case of
1-connected rational co-H-space.

On the other hand Hurvitz\cite{hurv} introduced the genus of maps
as follows

\begin{defn}
Let $f:X \to Y$ be a map between two nilpotent connected CW
  complexes $X,Y$. The  {\it genus}  $G(f)$ of $f$ is defined to be the
  set of equivalence classes of maps $f':X' \to Y'$ such that
  for each prime $p$ there exist homotopy equivalences $h_p:X'_{(p)} \to
  X_{(p)}$ and $k_p:Y'_{(p)} \to Y_{(p)}$ satisfying
   $f_{(p)}h_p \simeq k_pf'_{(p)}$, where two maps $f'_i :X'_i \to Y'_i$ , $i=1,2$, are {\it equivalent} if there
  exist homotopy equivalences $g:X'_1\to X'_2$ and $h:Y'_1 \to
  Y'_2$ such that $ f'_2g\simeq hf'_1$.
\end{defn}
Let $[f,f]_{t}$ be the set of equivalence classes of pairs $(h,k)$
of $t$-equivalences $h:X \to X$, $k:Y \to Y$ satisfying $kf \sim
fh$ where $t$ is a positive integer. One of the main results in
Hurvitz's paper is the following extension of Zabrodsky's result.
\begin{thm}
Let $f:X \to Y$ be an H-map between H-spaces such that
$H^*(X,\mathbb{Q})$ and $H^*(Y,\mathbb{Q})$ are primitively
generated and $H^*(f,\mathbb{Q})$ is either a monomorphism , an
epimorphism , an isomorphism or zero. Then $G(f)$ admits an
abelian group structure and there exist integers $k$ and $\hat{t}$
(depending on $X,Y$ and $f$) and an exact sequence
\[
[f,f]_{\hat{t}} \overset{\alpha'}{\to} [\mathbb{Z}^*_{\hat{t}}/\pm 1]^k
\overset{\hat{\xi}}{\to}  G(f) \to 0
\]
\end{thm}
The main result in this  paper is the following
\begin{thm}\label{T:main1}
Let $f:X \to Y$ be an H-map between H-spaces such that
$H^*(X,\mathbb{Q})$ and $H^*(Y,\mathbb{Q})$ are primitively
generated . Then $G(f)$ admits an abelian group structure and
there exist integers $k$ and $\hat{t}$ (depending on $X,Y$ and
$f$) and an exact sequence
\[
[f,f]_{\hat{t}} \overset{\alpha'}{\to} [\mathbb{Z}^*_{\hat{t}}/\pm 1]^k
\overset{\hat{\xi}}{\to}  G(f) \to 0
\]
where
\[ k= \left\{ \begin{array}{ll}
l(X) & \mbox{ if $H^*(f,\mathbb{Q})$ is an isomorphism }\\
l(X)+l(Y) & \mbox{ otherwise} \end{array} \right. \]
\newline and $l(X)$ is the number of $i$ such that $\pi_i(X)\bigotimes
\mathbb{Q}\neq 0$ while $\hat{t}$ will be defined in Section 2.
\end{thm}
\begin{rem}
Obviously the restriction on $H^*(f,\mathbb{Q})$ in Hurvitz's
result has been removed .
\end{rem}
It is not clear if Hurvitz's original approach can be extended to
get the dual result . Combining Hurvitz's and McGibbon's
approaches we are able to extend the above result to the dual case
as follows:
\begin{thm}\label{T:main2}
Let $f:X \to Y$ be a co-H-map between co-H-spaces with coprimitive
generated homotopy Lie algebras . Then $G(f)$ admits an abelian
group structure and there exist integers $k$ and $\hat{t}$
(depending on $X,Y$ and $f$) and an exact sequence
\[
[f,f]_{\hat{t}} \overset{\alpha'}{\to} [\mathbb{Z}^*_{\hat{t}}/\pm 1]^k
\overset{\hat{\xi}}{\to}  G(f) \to 0
\]
where
\[ k= \left\{ \begin{array}{ll}
l(X) & \mbox{ if $\pi_*(f)\bigotimes \mathbb{Q}$ is an isomorphism
}\\ l(X)+l(Y) & \mbox{ otherwise}
\end{array} \right.
 \]
 \newline and $l(X)$
is the number of $i$ such that $H_i(X,\mathbb{Q})\neq 0$ while
$\hat{t}$ will be defined in Section 3.
\end{thm}
\begin{rem}
There are similar results for $G_Y(f)$ and $G^X(f)$ as defined in
\cite{hurv}. We will omit these quite clear extensions.
\end{rem}
\begin{rem}
The analogue of rational (co-)H-space for map should be rational
(co-)H-map between rational (co-)H-spaces. Although we are unable
to extend the main results in this paper to the general case until
now, we do believe that it is true at least for rational
(co-)H-map between rational (co-)H-spaces with (co-)primitively
generated rational (homotopy Lie algebras) cohomology rings.
\end{rem}

All spaces considered in this paper are based and of the homotopy
type of 1-connected CW complexes of finite type. For simplicity
$0$ will be included into the set of primes. All spaces are
assumed to have finite number of nontrivial homotopy groups when
we are dealing with genus of H-map while spaces are assumed to
have finite number of nontrivial homology groups when we are
concerned with genus of co-H-map.

$\mathbb{Z}_t$ is the ring of $modt$ integers and $\mathbb{Z}^*_t$ is the group of
units in $\mathbb{Z}_t$. A map is said to be $t$-{\it equivalence} if it is $p$-equivalence for
all prime $p$ dividing $t$.

\section{Genus of maps between $H$-spaces}

First we fix some notations. Let $(X,\mu)$ be an H-space. We shall
denote by $+$ the operation on $[Z,X]$ induced by $\mu$ , by
$\phi_n$ the n-th power map
\[
\phi_n=\mu(\mu \times 1) \cdots (\mu \times 1 \times \cdots \times
1 ) \circ (\Delta \times 1 \times \cdots \times 1) \cdots \Delta
\]
where $\Delta$ denotes the diagonal map.

In general the product $\mu$  is not assumed to have an inverse.
Nevertheless we have the following elementary but important result
\begin{lem}
If $(X,\mu)$ is an H-space, then $[W,X]$ is an algebraic loop for
any space $W,$ i.e., for any two maps $f,g \in [W,X]$ there exists
a unique $D_{f,g} \in [W,X]$ such that $D_{f,g}+g=f$.
\end{lem}
\begin{lem}\label{T:lemma2}
Let $h:X_1 \to X_2$ be an H-map and $g:W_1 \to W_2$ be a map. For
any two maps $f_1$,$f_2 \in [W_2,X_1],$ we have
\[
hD_{f_1,f_2} = D_{hf_1,hf_2} \text{ and }
(D_{f_1,f_2})g=D_{f_1g,f_2g}.
\]
\end{lem}
\begin{rem}
For the proof of the lemmas above, see \cite{az1}
\end{rem}
For simplicity, denote $D_{f,g}$ by $f-g$. Define $$\bar{\mu}=\mu
-\pi_1 -\pi_2$$ where $\pi_1, \pi_2$ are projections from $X
\times X$ to its factors. Then for any coefficient $\mathbb{A}$ ,
define $$PH^*(X,\mathbb{A})=kerH^*(\bar{\mu},\mathbb{A})$$

 Now let $X$ be a rational H-space with
only finitely nontrivial homotopy groups. Denote by $K(X)$
the generalized Eilenberg-MacLane space with homotopy groups
$\pi_*(X)/torsion$. Let $\bar{\Delta}:X \to X \bigwedge X$ denote
the reduced diagonal map. Define $PH_*(X,\mathbb{Z})$ to be kernel
of $H_*(\bar{\Delta},\mathbb{Z})$ . Let $\sigma_n:\pi_n(X) \to
PH_*(X,\mathbb{Z})$ be the obvious quotient map of the Hurewicz
homomorphism. Since $X$ is a rational H-space , the map $\sigma_n$
has a finite kernel and cokernel for each $n$. Define
\[
 t_n(X)=exp(\text{coker}
 \sigma_{n+1})exp(\text{ker}\sigma_n)
\]
\[
t(X)=\prod_{n
\leq N(X)}t_n(X)
\]
where , for a given finite abelian group $G$, $exp{G}$ is the
smallest integer $n\geq 1$ such that $ng=0$ for all $g \in G$ and
$N(X)$ is the least integer such that $\pi_n(X)=0$ for every $n > N(X)$ .

 Let
$s_n(X)=exp(torsion(H^n(X,\mathbb{Z})))$ when
$QH^n(X,\mathbb{Q})\neq0$ and $s_n(X)=1$ otherwise. Let
$\mathbf{s}(X)$ be the sequence of integers $\{s_1(X), s_2(X)
\newline , \cdots\}$ . If $f:X \to Y$ is a map , then   $\mathbf{s}(f)$
is defined to be the sequence of integers $\{s_1(X)s_1(Y),
s_2(X)s_2(Y), \cdots\}$ .

On the other hand ,  given a   H-space $X$ with primitively
generated
 $H^*(X,\mathbb{Q})$ and a sequence of integers
  $\mathbf{s}=\{s_1, s_2, \cdots\}$ , an $\mathbf{s}$-maximal map
for $X$ is a map $\varphi_X:X \to K(X)$ such that for each $n$ there exist
 $\{x_1, \cdots , x_r\} \in H^n(X,\mathbb{Z})$ which projects to a basis
  of $PH^n( X,\mathbb{Q})$ and a basis $\{\iota_1, \cdots , \iota_r\}$
   for $PH^n( K(X),\mathbb{Z})/torsion$ with $(\varphi_X)^*(\iota_i)=s_nx_i$.
    If $X$ is an H-space and  $\varphi_X$ is an H-map, then the maximal
     map will be called primitive. An elementary but  useful result is
as follows:
\begin{lem}
   $s_nt_n(X)$ can be divided by the integer
$exp(\pi_n(fiber(\varphi_X)))$ for  an $\mathbf{s}$-maximal map of
$X$ .
\end{lem}
\begin{rem}
Any $\mathbf{s}(X)$-maximal map of an H-space $X$ is  primitive.
\end{rem}
Now we define the integer $\hat{t}$  in
Theorem\ref{T:main1} as follow:
\begin{defn}
Given an  H-map $f:X \to Y$ between two H-spaces,
 $$ \hat{t}=t(X)(t(Y))^2\underset{n \leq \max(N(X),N(Y))}
 {\prod}s_n(X)s_n(Y) $$
\end{defn}

 Let $f:X\to Y$ be a map between two rational H-spaces
 with primitively generated rational cohomology and $\mathbf{s}$ be a
 sequence of integers $\{s_1, s_2, \cdots\}$. Then an
 $\mathbf{s}$-maximal map for $f$ is a homotopy commutative diagram
 \[
 \begin{CD}
 X @>f>> Y \\ @V\varphi_XVV @V\varphi_YVV \\
 K(X) @>f_X>> K(Y)
 \end{CD}
 \]
 where $\varphi_X, \varphi_Y$ are $\mathbf{s}$-maximal maps. If the map
  $f$ is also an H-map, then an  $\mathbf{s}$-maximal map for
   $f$ is called primitive if $f_X, \varphi_X, \varphi_Y$ are H-maps.
\begin{rem}
Any  $\mathbf{s}(f)$-maximal map for $f$ is primitive.
\end{rem}

In this section we will prove Theorem \ref{T:main1}. Before that ,
however, we need several preliminary results .
\begin{thm}\label{T:finite-to1-H}
Let $X$ be an H-space with primitively generated
$H^*(X,\mathbb{Q})$ and with only finitely nontrivial homotopy
groups and $W$ be any space. Then
\newline (i)For every map $f:W_p \to X_p$, there exist an integer
$n$ and a map $g:W \to X$ such that $(n,p)=1$ and $g_{(p)} \sim
\phi_n f$.
\newline (ii)Given two maps $f_1,f_2:W \to X$ so that $\varphi_{X}f_1 \sim \varphi_X f_2$
where $\varphi_{X}:X \to K(X)$ is a primitive
$\mathbf{s}$-maximal map, then $\phi_{t'} f_2 \sim \phi_{t'}f_1$
where $t'=\underset{n \leq N(X)}{\prod}t_n(X)s_n$.
\end{thm}
\begin{proof}
The statement (i) is essentially the same as that of Theorem 2.1
in \cite{hurv}.

 To prove the part (ii),it suffices to prove that $\phi_{t'} k \sim *$ by
Lemma \ref{T:lemma2} where $k=D_{f_1,f_2}$. The given conditions
and Lemma \ref{T:lemma2} imply that $ \varphi_{X}k\sim *$. The
result follows from an induction by Postnikov tower of map $
\varphi_X$ , see the proof of Theorem \ref{T:finite-to1-co-H} for
details of a dual proof.
\end{proof}
\begin{thm}\label{T:2.2}
Let $f:X \to Y$ be an H-map between two H-spaces and $l$ be an
integer with $l=p_1^{w_1}\cdots p_s^{w_s}$ . Than given two spaces
$X',Y'$, a function $f':X' \to Y'$ and homotopy equivalences
$h_{p_i}:X'_{(p_i)} \to X_{(p_i)}$, $k_{p_i}:Y'_{(p_i)} \to
Y_{(p_i)}$ satisfying $f_{(p_i)}h_{p_i} \sim k_{p_i}f'_{(p_i)}$,
 there exist $l$-equivalences $h:X' \to X$, $Y' \to Y$ so
that $fh \sim kf'$.
\end{thm}
\begin{proof}
It is Theorem 2.2 in \cite{hurv}
\end{proof}
\begin{thm}\label{T:3.4}
If $f:X \to Y$ is any map between rational H-spaces. Suppose that
an $\mathbf{s}(f)$-maximal map of $f$ is given as in the following
diagram.
\[
\begin{diagram}
 \node[2]{X'} \arrow[2]{e,t}{h} \arrow{s,r}{\varphi_{X'}}
\arrow{ssw,l}{f'} \node[2]{X} \arrow{s,r}{\varphi_X}
\arrow{ssw,r}{f}\\ \node[2]{K(X)} \arrow[2]{e,t}{\alpha}
\arrow{ssw,l}{g} \node[2]{K(X)} \arrow{ssw,l}{g}\\ \node{Y'}
\arrow{s,l}{\varphi_{Y'}} \arrow[2]{e,t}{k} \node[2]{Y}
\arrow{s,l}{\varphi_Y}\\ \node{K(Y)} \arrow[2]{e,t}{\beta}
\node[2]{K(Y)}
\end{diagram}
\]
Then for every pair $(\alpha,\beta)$, where $\alpha:K(X) \to K(X)$
and $\beta:K(Y) \to K(Y) $ are $\hat{t}$-equivalences such that
$\beta g \sim g\alpha$, there exist a map $f':X' \to Y' , f' \in
G(f)$ and $\hat{t}$-equivalences $h:X' \to X$ and $k:Y' \to Y$
such that the every diagram in the following cube commutes up to
homotopy and $\varphi_{X'},\varphi_{Y'}$ are
$\mathbf{s}(f)$-maximal maps.
\end{thm}
\begin{rem}\label{T:3.4re}
It is easy to know that, if $\alpha,\beta$ are homotopy equivalences,
then $f'$ is equivalent to $f$.
\end{rem}
\begin{proof}
It is Theorem 3.4 in \cite{hurv}.
\end{proof}
\begin{thm}\label{T:realiz}
Let $f:X_1 \to X_2$ be an H-map between H-spaces such that
$H^*(X_1,\mathbf{Q})$ and $H^*(X_2,\mathbf{Q})$ are primitively
generated. Let $f':Y_1 \to Y_2$ be a map between rational
H-spaces. Choose $\mathbf{s}(f)$-maximal maps $f_X , f'_Y $ for
$f, f'$ respectively so that the maximal map for $f$ is primitive.

Then given a pair of maps $(\alpha_1,\alpha_2)$ as in the
following diagram satisfying $\alpha_2 f'_Y = f_X \alpha_1$, there
exist mappings $h_i:Y_i \to X_i$, $i=1,2$, so that all the
diagrams in the following cube commute up to homotopy.
\[
\begin{diagram}
 \node[2]{Y_1} \arrow[2]{e,t}{h_1} \arrow{s,r}{\varphi_{Y_1}}
\arrow{ssw,l}{f'} \node[2]{X_1} \arrow{s,r}{\varphi_{X_1}}
\arrow{ssw,r}{f}\\ \node[2]{K(Y_1)}
\arrow[2]{e,t}{\phi_{\hat{t}}\alpha_1} \arrow{ssw,l}{f'_Y}
\node[2]{K(X_1)} \arrow{ssw,l}{f_X}\\ \node{Y_2}
\arrow{s,l}{\varphi_{Y_2}} \arrow[2]{e,t}{h_2} \node[2]{X_2}
\arrow{s,l}{\varphi_{X_2}}\\ \node{K(Y_2)}
\arrow[2]{e,t}{\phi_{\hat{t}}\alpha_2} \node[2]{K(X_2)}
\end{diagram}
\]
\end{thm}
\begin{proof}
By Proposition 1.8 in \cite{az2} , for each pair
$(\alpha_1,\alpha_2)$ in the Theorem, there exist mappings
$h'_i:Y_i \to X_i$ so that $\varphi_{X_i}h'_i \sim
\varphi_{Y_i}\phi_{t(X_1)t(X_2)}\alpha_i$. It is easy to know that
$\varphi_{X_2} h'_2 g \sim \varphi_{X_2} f h'_1$ . Theorem
\ref{T:finite-to1-H} implies that $$\phi_{(t(X_2))^2} h'_2 g \sim
\phi_{(t(X_2))^2} f h'_1\sim f\phi_{(t(X_2))^2}h'_1$$ since $f$ is an
H-map. The result follows by taking $h_i =\phi_{t'}h'_i $  where
$t'=\hat{t}/t(X_1)t(X_2)$.
\end{proof}
\begin{thm}\label{T:4.6}
Let $f$ be as in Theorem \ref{T:main1}. Given a primitive
$\mathbf{s}(f)$-maximal map $g$ of $f$, there exists a commutative
cube up to homotopy  as follows
\[
\begin{diagram}
 \node[2]{X\times K(X)} \arrow[2]{e,t}{\vartheta_X} \arrow{s,r}{\varphi_{X}\times \phi_{\hat{t}}}
\arrow{ssw,l}{f \times g} \node[2]{X} \arrow{s,r}{\varphi_X}
\arrow{ssw,r}{f}\\ \node[2]{K(X)\times K(X)} \arrow[2]{e,t}{\mu}
\arrow{ssw,l}{g\times g} \node[2]{K(X)} \arrow{ssw,l}{g}\\
\node{Y\times K(Y)} \arrow{s,l}{\varphi_{Y}\times \phi_{\hat{t}}}
\arrow[2]{e,t}{\vartheta_Y} \node[2]{Y} \arrow{s,l}{\varphi_Y}\\
\node{K(Y)\times K(Y)} \arrow[2]{e,t}{\mu} \node[2]{K(Y)}
\end{diagram}
\]
where $\vartheta_X, \vartheta_Y$ restrict to the identity on the
first factor and $\mu$ are the standard multiplication on product
of Eilenberg-MacLane spaces and left face is the given maximal
map.
\end{thm}
\begin{proof}
The proof is similar to that of Theorem \ref{T:realiz} using
Proposition 4.6 in \cite{mcg}. Although we  modified the notion of
maximal map,  we have also modified the integer $\hat{t}(X)$ to
ensure that it remains valid.
\end{proof}
\begin{thm}\label{T:modif}
Let $f:X \to Y$ be as in Theorem \ref{T:main1} . Given $[f'] \in
G(f)$, then there is a commutative cube up to homotopy
\[
\begin{diagram}
 \node[2]{X'} \arrow[2]{e,t}{h} \arrow{s,r}{\varphi_{X'}}
\arrow{ssw,l}{f'} \node[2]{X} \arrow{s,r}{\varphi_X}
\arrow{ssw,r}{f}\\ \node[2]{K(X)} \arrow[2]{e,t}{\alpha}
\arrow{ssw,l}{g} \node[2]{K(X)} \arrow{ssw,l}{g}\\ \node{Y'}
\arrow{s,l}{\varphi_{Y'}} \arrow[2]{e,t}{k} \node[2]{Y}
\arrow{s,l}{\varphi_Y}\\ \node{K(Y)} \arrow[2]{e,t}{\beta}
\node[2]{K(Y)}
\end{diagram}
\]
where $h,k$ are $\hat{t}$-equivalences , all the vertical maps are
maximal maps, the back and front face diagrams are homotopy
pullback diagrams and
$$\alpha=\prod_i\alpha_i:K(X)\to K(X),\beta=\prod_j\beta_j:K(Y)\to K(Y)$$ where
$$\alpha_i:K(\pi_i(X)/torsion,i)\to
K(\pi_i(X)/torsion,i)$$$$\beta_j:K(\pi_j(Y)/torsion,j)\to
K(\pi_j(Y)/torsion,j)$$ and
$$K(X)=\prod_iK(\pi_i(X)/torsion,i),K(Y)=\prod_jK(\pi_j(Y)/torsion,j)$$
\end{thm}
\begin{proof}
Take $\hat{t}$-equivalences $h:X' \to X, k:Y' \to Y$ as granted by
Theorem \ref{T:2.2} so that $fh \sim kf'$.  We will try to define
the other maps in the above diagram so that the Theorem is true.
Without loss of generality, we can assume that there is only one
$l$ such that $\pi_i(X)\bigotimes \mathbb{Q} =\pi_i(Y)\bigotimes
\mathbb{Q} =0$ if $i\neq l$. Choose rational H-space structures on
$X',Y'$ such that $h_{(0)},k_{(0)},f'_{(0)}$ are H-maps. It is
always possible to choose bases$\{x_1,\cdots,x_{r_1}\}$ for
$PH^l(X,\mathbb{Z})/torsion$ and $\{x'_1,\cdots,x'_{r_2}\}$ for
$PH^l(Y,\mathbb{Z})/torsion$  so that $PH^*(f,\mathbb{Z})/torsion$
is diagonal with respect to these bases and
$\{x_1,\cdots,x_{r_1}\}$ and $\{x'_1,\cdots,x'_{r_2}\}$ are  bases
for $PH^*(X,\mathbb{Q})$ and $PH^*(Y,\mathbb{Q})$.  Let
$\varphi_X,\varphi_Y$ be defined by
$s_l(f)x_i=(\varphi_X)^*(\iota_i)$ and
$s_l(f)x'_j=(\varphi_Y)^*(\iota'_j)$  where
$\{\iota_1,\cdots,\iota_{r_1}\}$ is a basis for
$H^l(K(X),\mathbb{Z})/torsion$ and
$\{\iota'_1,\cdots,\iota'_{r_2}\}$ is a basis for
$H^l(K(Y),\mathbb{Z})/torsion$. Then
 $g$ is also of the diagonal form under the bases. Choose elements
 $\{y_1,\cdots,y_{r_1}\} \in H^l(X',\mathbb{Z})$ and $\{y'_1,\cdots,y'_{r_2}\} \in H^l(Y',\mathbb{Z})$ such that they projects to
 bases of $PH^*(X',\mathbb{Q})$ and
$PH^*(Y',\mathbb{Q})$ respectively. Define the maps $\varphi_{X'}, \varphi_{Y'}$ similarly. Since $k,h$ are
$\hat{t}$-equivalences and rational H-maps, we can write
\[
(h)^*(x_i)=\sum\lambda_{ii'} y_{i'} + v_i ,
(k)^*(x'_j)=\sum\lambda'_{jj'} y'_{j'} +w_j
\]
where $\det\lambda_{ii'},\det\lambda'_{jj'}$ are prime relative to
$\hat{t}$ and $v_i,w_j$ are torsions. Certainly we can find
$\bar{v}_i$ and $\bar{w}_j$ so that
\[
(h)^*(x_i)=\sum\lambda_{ii'} (y_{i'} + \bar{v}_{i'} ),
(k)^*(x'_j)=\sum\lambda'_{jj'} (y'_{j'} +\bar{w}_{j'})
\]

Let $\varphi_{X'}$ be the map which sends each $\iota_i$ to
$s_l(f)(y_i +\bar{v}_i)$ and $\varphi_{Y'}$  be the map which
sends each $\iota'_j$ to $s_l(f)(y'_j +\bar{w}_j)$ and set
$\beta,\alpha$ to be the maps which send $\iota_i$ and $\iota'_j$
to $\sum\lambda_{ii'} \iota_{i'}$ and
 $\sum\lambda'_{jj'} \iota'_{j'}$ respectively. Obviously what
we have to verify is that $\varphi_{Y'} f' \sim g\varphi_{X'} $.
It is easy to know from diagram chasing argument that $\beta
\varphi_{Y'} f' \sim \beta g\varphi_{X'} $. Now $(\varphi_{Y'}
f')^*(\iota_i)-(g\varphi_{X'})^*(\iota_i)=u_i$ which are torsion
of order dividing $\hat{t}$ for all $i$. $u_i=0$ for all $i$ iff
$\sum \lambda_{ii'}u_{i'}=0$ for all $i$ iff $\beta \varphi_{Y'}
f' \sim \beta g\varphi_{X'}$ since $\det\lambda_{ii'}$ is prime
relative to $\hat{t}$. Now it is clear that the front and back
squares are homotopy pullback.
\end{proof}
\begin{proof}[Proof of Theorem \ref{T:main1}]
Without loss of generality, we can assume that there is an integer
$l$ such that $\pi_i(X)\bigotimes \mathbb{Q}=\pi_i(Y)\bigotimes \mathbb{Q}=0$ if
$i\neq l$. By Theorem \ref{T:3.4} and Theorem \ref{T:modif} it
follows that there exists a surjection $\xi:T \to G(f)$ where
\[
T=\{(\alpha,\beta)|\beta g \sim g\alpha , \text{$\alpha$ and
$\beta$ are $\hat{t}$-equivalences}\} .
\]
In other words, there is a surjection $\xi':T' \to G(f)$ where
\[
T'=\{(A_1,A_2)|A_i \in M(\mathbb{Z},r_i),(\det A_i,\hat{t} )=1,
i=1,2, A_2C=CA_1 \}
\]
where $M(\mathbb{Z},n)$ is the set of $n\times n$ matrices and
$A_1,A_2,C$ represents $\alpha,\beta,g$ respectively.

Given $(A_1,A_2), (B_1,B_2) \in T^{'}$ such that $\det A_i \equiv
\det B_i $ (mod $\hat{t}$) for $i=1,2$, then an elementary
calculation shows that $A_1 =B_1 G_1+\hat{t}H_1$, $A_2=B_2 G_2
+\hat{t}H_2$, $(G_1,G_2)\in T^{'}, \det G_1 \equiv \det G_2 \equiv
1(mod\hat{t})$ and $H_2C=CH_1$. Claim \ref{T:claim} implies that
we can assume $\det G_1=\det G_2=1$. It follows from Theorem
\ref{T:modif} and Theorem \ref{T:4.6} and Remark \ref{T:3.4re}
that $\xi^{'}(A_1,A_2)=\xi^{'}(B_1,B_2)$. Thus $\xi^{'}$ factors
through the map $\det:T^{'} \to
(\mathbb{Z}^*_{\hat{t}})^{l(X)+l(Y)}$. The rest of the proof is
the same as that of Hurvitz\cite{hurv}.
\end{proof}

\begin{claim}\label{T:claim}
Given $G_1, G_2$ with $\det G_1 \equiv \det G_2 \equiv 1
(mod\hat{t})$ and $(G_1,G_2) \in T'$. There exist $H_1 \in
GL(\mathbb{Z},r_1)$, $H_2 \in GL(\mathbb{Z},r_2)$ such that
$(H_1,H_2) \in T'$ and $G_1 H_1 = Id + \hat{t}H'_1$, $ H_2 G_2 =
Id + \hat{t}H'_2$
\end{claim}
\begin{proof}
Assume, without loss of generality, that the $m_1\times m_2$ matrix
$C$ has
 entries zero unless those on diagonal.  Then $(G_1,G_2) \in T'$ is
  equivalent to the equations
\[
(G_1)_{ij}c_{jj}=c_{ii}(G_2)_{ij} , \text{  for  } 0 \leq i,j \leq
l_0
\]
\[
(G_2)_{ij}=0  \text{  for  } i\leq l_0 \text{  and  } j>l_0
\]
\[
(G_1)_{ij}=0  \text{  for  } j\leq l_0 \text{  and  } i>l_0
\]
where $c_{ii}$,  $0 \leq i \leq l_0$, are the nontrivial terms in
matrix $C$. The proof of the existence of the matrices $H_1, H_2$
with the prescribed properties will be given by mathematical
induction and it is obvious for $l=1$. Assume it is also true for
$l < s$, we will prove it for $l=s$. If $l_0 < m_1$, using
elementary matrix we can find
 $H_1 \in GL(\mathbb{Z},r_1)$  such that  $(H_1,Id) \in T'$ and
 we have the following modulo  $\hat{t}$ equation
\[
(G_1 H_1)_{ij}= \left\{ \begin{array}{ll}
 0  & \mbox{ if $i=m_1 , l_0 <j<m_1$ }\\
 1  & \mbox{ if $i=m_1 , j=m_2$ }\\
 (G_1)_{ij} & \mbox{ otherwise}
 \end{array}
 \right.
\]
 Similarly we can find $H_2 \in GL(\mathbb{Z},r_2)$ if $l_0 <
 m_2$ with similar properties.
It follows that we can assume $l_0 = m_1 =m_2$. In this case
$(G_1,G_2) \in T'$ is equivalent to the equations
\[
(G_1)_{ij}c_{jj}=c_{ii}(G_2)_{ij} , \text{  for  }  i,j \leq l_0
\]
The equations above imply that $G_1$ and $G_2$ are determined by one
matrix. Thus the result follows from the corresponding result in the
one matrix case.
\end{proof}

\section{Genus of co-H-maps between co-H-spaces}
As in the last section we  fix some notations first. Throughout
this section all spaces will be 1-connected finite CW complexes
with rational co-H-space structures. Given space $X$, let $M(X)$
denote the bouquet of spheres whose reduced integral homology is
isomorphic to $\tilde{H}_*(X,\mathbb{Z})/torsion$. For any space
$X$, let $Q\pi_*(X)$ denote the quotient  of $\pi_*(X)$ by the
subgroup generated by all Whitehead products in $\pi_*(X)$. An
$\mathbf{s}$-maximal map from $M(X)$ to $X$ is one that satisfies
a condition dual to that in the rational H-space case. A
$\mathbf{s}$-maximal map of a co-H-space is called coprimitive if
it is a co-H-map.

Let  $s_n(X)=exp(torsion(\pi_k(X)))$ when
 $Q\pi_n(X)\bigotimes \mathbb{Q}\neq0$ and $s_n(X)=1$ otherwise.
  Let $\mathbf{s}(X)$ be the sequence of integers
   $\{s_1(X), s_2(X), \cdots\}$ . If $f:X \to Y$ is a map ,
   then   $\mathbf{s}(f)$ is defined to be the sequence of integers
    $\{s_1(X)s_1(Y), s_2(X)s_2(Y), \cdots\}$ .

 On the other hand , let $\sigma_n:
Q\pi_n(X)/torsion \to H_n(X,\mathbb{Z})$ be induced by the
Hurewicz homomorphism. Since $X$ is a rational co-H-space,
$\sigma_n$ has finite kernel and finite cokernel for each $n$.
Define
\[
 t_n(X)=exp(\text{coker}
 \sigma_{n+1})exp(\text{ker}
 \sigma_{n})
\]
and
\[
t(X)=\prod_{n
\leq N(X)}t_n(X)
\]
where $N(X)$ is the least integer so that , for every $n >
N(X), H_n(X)=0$. As in the last section, we have the following
\begin{lem}
$s_nt_n(X)$ can be divided by the integer
 $exp(H_n(C(\varphi_X),\mathbb{Z}))$
where $C(f)$ is the homotopy cofiber of $f$ and
 $\varphi_X$ is an $\mathbf{s}$-maximal map of $X$.
\end{lem}
\begin{rem}
Any $\mathbf{s}(X)$-maximal map of a co-H-space is also
coprimitive.
\end{rem}
 Now we define as before an integer $\hat{t}$ in
Theorem \ref{T:main2} as follow
\begin{defn}
Given a  co-H-map $f:X \to Y$ between two co-H-spaces, $$
\hat{t}=t(X)(t(Y))^2\underset{n \leq \max{(N(X),N(Y))}}
{\prod}s_n(X)s_n(Y) $$
\end{defn}
Let $X$ be a rational co-H-space, we say $\pi_*(X)\bigotimes Q$ is
coprimitive generated as Lie algebra  if there exists a Lie
algebra basis $\{x_1,\cdots,x_n\} \in \pi_*(X)\bigotimes Q \cong
\pi_*(X_{(0)})$ such that $x_i$ is represented by a co-H-map $S^j
\to X_{(0)}$.

 Let $f:X\to Y$ be a map between two rational co-H-spaces with
  coprimitive generated rational homotopy groups. Then an
  $\mathbf{s}$-maximal map for $f$ is a homotopy commutative diagram
 \[
 \begin{CD}
 X @>f>> Y \\ @A\varphi_XAA @A\varphi_YAA \\
 M(X) @>f_X>> M(Y)
 \end{CD}
 \]
 where $\varphi_X, \varphi_Y$ are $\mathbf{s}$-maximal maps.
 If the map $f$ is also a co-H-map, then an  $\mathbf{s}$-maximal map for  $f$ is called coprimitive
 if $f_X, \varphi_X, \varphi_Y$ are co-H-maps.
\begin{rem}
Any  $\mathbf{s}(f)$-maximal map of a co-H-map $f$ is coprimitive.
\end{rem}

In this section we will prove Theorem \ref{T:main2}. As before , we
need several preliminary results which are dual to those in the
previous section.
 Let $(X,\nu)$ be a
co-H-space. We shall denote by $+$ the operation on $[X,W]$
induced by $\nu$ , by $\eta_n$ the  map
\[
\eta_n=F \cdots (F \vee 1 \vee \cdots \vee 1 ) \circ (\psi \vee 1
\vee \cdots \vee 1) \cdots(\psi \vee 1) \psi
\]
where F denotes the folding map.

In general the product $\nu$  is not assumed to have an inverse.
Nevertheless we have the following elementary but important result
\begin{lem}
If $(X,\nu)$ is a co-H-space, then $[X,W]$ is an algebraic loop
for any space $W$, i.e., for any pair $f,g \in [X,W]$ there exists
a unique $D_{f,g} \in [X,W]$ such that $D_{f,g}+g=f.$
\end{lem}
\begin{lem}\label{T:lemma32}
Let $h:X_1 \to X_2$ be a co-H-map and $g:W_1 \to W_2$ be a map.
For any pair $f_1,f_2 \in [X_2,W_1],$ we have $(D_{f_1,f_2})h =
D_{f_1h,f_2h}$ and $gD_{f_1,f_2}=D_{gf_1,gf_2}.$
\end{lem}

\begin{rem}
The proof of the lemmas above are dual to that in  \cite{az1}
\end{rem}
The dual of Theorem \ref{T:finite-to1-H} is
\begin{thm}\label{T:finite-to1-co-H}
Let $X$ be a 1-connected co-H-space and $W$ be a space with a
finite number of nontrivial homotopy groups. Then we have the following: \newline
(i)For every map $f:X_{(p)} \to W_{(p)}$ there exist an integer
$n$, $(n,p)=1$ , and a map $\bar{f}:X \to W$ so that
$\bar{f}_{(p)} \sim f \eta_n $.\newline(ii) If two maps
$f_1,f_2:X \to W$ satisfy $f_1 \varphi_{X} \sim f_2 \varphi_{X}$,
 then $ f_2 \eta_{t'} \sim f_1 \eta_{t'}$  where $\varphi_{X}:M(X) \to X$ is a coprimitive
$\mathbf{s}$-maximal map and $t'=\underset{n \leq
N(X)}{\prod}t_n(X)s_n$.
\end{thm}
The part (i) is dual to  Theorem 2.1(i) in \cite{hurv}, so we will
only give a detailed proof of the part (ii) .
\begin{proof}

As mentioned before, the maximal map $\varphi_X$ as in the Theorem
exists. Let $k=D_{f_1,f_2}$. It suffices to prove that $ k
\eta_{t'} \sim *$ by Lemma \ref{T:lemma32}.  The given condition and Lemma \ref{T:lemma32} imply that $k \varphi_{X}\sim
*$. Now consider the
 Moore-Postnikov
decomposition of the map $\varphi_{X}:M(X) \to X$. A typical term
in the decomposition fits into the following
\[
\begin{CD}
X^{(n)} @>g^{(n)}>> X^{(n+1)} \\ @ AAA @VVV \\ M(X)
@>\varphi_X>>  X
\end{CD}
\]
where $X^{(0)}= M(X)$ and map $X^{(n)} \to X^{(n+1)} $ is the
mapping cone of a map $k^{(n)}: M(H_{n+1}(C(\varphi_{X}),\mathbb{Z}),n) \to X^{(n)}$
and $M(G,n)$ is the type $G-n$ Moore space :
\[
\tilde{H}_i(M(G,n),\mathbb{Z})= \left\{ \begin{array}{ll}
                              G  & \mbox{if $i=n$}\\
                              0  & \mbox{if $i\neq n$}
                              \end{array}
                       \right.
                           \]
By the condition on $X$ , there is an integer $s$ such that
$X^{(s)}=X$.  Let $t_n=\prod_{i \leq n}t_i(X)s_i$. We will prove
that $k_n \eta_{t_n} \sim *$ where $k_n : X^{(n)} \to W$ is the
map induced by $k:X \to W$ which will completes the proof since
for $n=s$, $k=k_n$ and $t'=t_n$. We will proceed by induction. It
is obvious  for the case $n=0$ . If the statement is true for $n$
, we want to prove it for $n+1$. From the cofibration $$
M(H_{n+1}(C(\varphi_X),\mathbb{Z}),n) \to X^{(n)} \to X^{(n+1)}$$
the following exact sequence is exact:
\[
[M(H_{n+1}(C(\varphi_X),\mathbb{Z}),n+1),W] \to [ X^{(n+1)},W] \to
[X^{(n)},W]
\]
Now $(\eta_{t_n})^*(g^{(n)})^*(k_{n+1})=(\eta_{t_n})^* k_n =0$ .
Since $g^{(n)}$ is a co-H-map, it follows by Theorem \ref{T:co-H}
that $k_{n+1}\eta_{t_n}$ belongs to image of the group \newline
$[M(H_{n+1}(C(\varphi_X),\mathbb{Z}),n+1),W]$ which has exponent
dividing $t_{n+1}(X)$. Therefore it follows  that
$k_{n+1}\eta_{t_{n+1}s_{n+1}} \sim *$ which completes the
induction.
\end{proof}
The result dual to Theorem \ref{T:2.2} is
\begin{thm}\label{T:dual2.2}
Let $f:X \to Y$ be a co-H-map between two co-H-spaces, $l$ be an
integer with $l=p_1^{w_1}\cdots p_s^{w_s}$ . Given two spaces
$X',Y'$, a function $f':X' \to Y'$ and homotopy equivalences
$h_{p_i}:X_{p_i} \to X'_{p_i}$, $k_{p_i}:Y_{p_i} \to Y'_{p_i}$
satisfying $ k_{p_i}f_{(p_i)} \sim f'_{(p_i)}h_{p_i} $, then there
exist $l$-equivalences $h:X \to X'$, $Y \to Y'$ so that $kf \sim
f'h$.
\end{thm}
\begin{proof}
The proof is exactly the dual of that of Theorem 2.2 in
\cite{hurv}.
\end{proof}
The following is what is dual to Theorem \ref{T:3.4}
\begin{thm}\label{T:dual3.4}
If $f:X \to Y$ is any map  between rational co-H-spaces , Suppose
that an $\mathbf{s}(f)$-maximal map of $f$ is given as in the
following diagram.
\[
\begin{diagram}
 \node[2]{M(X)} \arrow[2]{e,t}{\alpha} \arrow{s,r}{\varphi_{X}}
\arrow{ssw,l}{g} \node[2]{M(X)} \arrow{s,r}{\varphi_{X'}}
\arrow{ssw,r}{g}\\ \node[2]{X} \arrow[2]{e,t}{h} \arrow{ssw,l}{f}
\node[2]{X'} \arrow{ssw,l}{f'}\\ \node{M(Y)}
\arrow{s,l}{\varphi_{Y}} \arrow[2]{e,t}{\beta} \node[2]{M(Y)}
\arrow{s,l}{\varphi_{Y'}}\\ \node{Y} \arrow[2]{e,t}{k}
\node[2]{Y'}
\end{diagram}
\]

Then for every pair $(\alpha,\beta)$ where
$\alpha:M(X) \to M(X)$, $\beta:M(Y) \to M(Y) $ are
$\hat{t}$-equivalences  such that $\beta g \sim g\alpha$, there
exist a map $f':X' \to Y' , f' \in G(f)$ and
$\hat{t}$-equivalences $h:X \to X'$ and $k:Y \to Y'$ so that the
every diagram in the above cube commutes up to homotopy and
$\varphi_{X'},\varphi_{Y'}$ are $\mathbf{s}(f)$-maximal maps.
\end{thm}
\begin{rem}\label{T:3.4re-co}
It is easy to know that, if $\alpha, \beta$ are homotopy equivalences,
then $f'$ is equivalent to $f$.
\end{rem}
\begin{proof}
The proof is dual to that of Theorem 3.4 in \cite{hurv}.
\end{proof}
That Dual to Theorem \ref{T:realiz} is the following
\begin{thm}\label{T:realiz3}
Let $f:X_1 \to X_2$ be a co-H-map between co-H-spaces such that
$\pi_*(X_1)\bigotimes \mathbb{Q}$ and $\pi_*(X_2)\bigotimes
\mathbb{Q}$ are coprimitively generated.  Let $f':Y_1 \to Y_2$ be
a map
  between rational co-H-spaces.
  Choose $\mathbf{s}(f)$-maximal maps $f_X , f'_Y $ for $f, f'$,
respectively, so that the maximal map for $f$ is  coprimitive.
Then given a pair of maps $(\alpha_1,\alpha_2)$ satisfying
$\alpha_2 f_X = f'_Y \alpha_1$ as in the following diagram, there
exist mappings $h_i:X_i \to Y_i$, $i=1,2$, so that all the
diagrams in the following cube commute up to homotopy.
\[
\begin{diagram}
 \node[2]{M(X_1)} \arrow[2]{e,t}{\alpha_1 \eta_{\hat{t}}} \arrow{s,r}{\varphi_{X_1}}
\arrow{ssw,l}{f_X} \node[2]{M(Y_1)} \arrow{s,r}{\varphi_{Y_1}}
\arrow{ssw,r}{f'_Y}\\ \node[2]{X_1} \arrow[2]{e,t}{h_1}
\arrow{ssw,l}{f} \node[2]{Y_1} \arrow{ssw,l}{f'}\\ \node{M(X_2)}
\arrow{s,l}{\varphi_{X_2}} \arrow[2]{e,t}{\alpha_2 \eta_{\hat{t}}}
\node[2]{M(Y_2)} \arrow{s,l}{\varphi_{Y_2}}\\ \node{X_2}
\arrow[2]{e,t}{h_2} \node[2]{Y_2}
\end{diagram}
\]
\end{thm}
\begin{proof}
An exactly dual proof can be given for Theorem \ref{T:realiz} ,
except that,  instead of  using Theorem \ref{T:finite-to1-H} and Proposition
1.8 in \cite{az2} , we appeal to   Theorem \ref{T:finite-to1-co-H}
and Theorem \ref{T:dual1.8} below which is dual to Proposition 1.8
in \cite{az2}.
\end{proof}
\begin{thm}\label{T:dual1.8}
Let $X$ be a finite co-H-space  and $\varphi_X:M(X) \to X$ be a
$\mathbf{s}$-maximal map and $g:M(X) \to W$ be any map into a space
$W$. Then there exists a map $h:X \to W$ such that $h \varphi_X
\sim g \eta_{t(X)}$.
\end{thm}
\begin{proof}
Let $\vartheta:X \to X \bigvee M(X)$ be the co-action given in
Proposition 5.6 \cite{mcg} with respect to maximal map
$\varphi_X$. Then $h$ is the composite $X \overset{\vartheta}{\to}
X\bigvee M(X)\overset{g\eta_{t(X)}\bigvee
\text{id}}{\longrightarrow}W\bigvee X \overset{\nabla}{\to}W$
where $\nabla|_{W}= \text{id}$ and $\nabla|_{X} = *$. It is easy to
verify that $h$ is what we want.
\end{proof}
\begin{thm}\label{T:4.6-co}
Let $f$ be as in Theorem \ref{T:main2}. Given a coprimitive
$\mathbf{s}(f)$-maximal map $g$ of $f$, there exists a commutative
cube up to homotopy  as follows
\[
\begin{diagram}
 \node[2]{M(X)} \arrow[2]{e,t}{\nu} \arrow{s,r}{\varphi_{X}}
\arrow{ssw,l}{g} \node[2]{M(X)\bigvee M(X)}
\arrow{s,r}{\varphi_X\bigvee \eta_{\hat{t}}}
\arrow{ssw,l}{g\bigvee g}\\ \node[2]{X}
\arrow[2]{e,t}{\vartheta_X} \arrow{ssw,l}{f} \node[2]{X\bigvee
M(X)} \arrow{ssw,r}{f\bigvee g}\\ \node{M(Y)}
\arrow{s,l}{\varphi_{Y}} \arrow[2]{e,t}{\nu} \node[2]{M(Y)\bigvee
M(Y)} \arrow{s,l}{\varphi_Y\bigvee \eta_{\hat{t}}}\\ \node{Y}
\arrow[2]{e,t}{\vartheta_X} \node[2]{Y\bigvee M(Y)}
\end{diagram}
\]
where $\vartheta_X, \vartheta_Y$ project to the identity on the
first factor and $\nu$ are the standard multiplication on bouquet
of spheres.
\end{thm}
\begin{proof}
The proof is similar to that of Theorem \ref{T:realiz} using
Proposition 5.6 in\cite{mcg}.
\end{proof}

\begin{thm}\label{T:modifco}
Let $f:X \to Y$ be as in Theorem \ref{T:main2} . Given  $[f'] \in
G(f)$, then there is a commutative cube
 up to homotopy
\[
\begin{diagram}
 \node[2]{M(X)} \arrow[2]{e,t}{\alpha} \arrow{s,r}{\varphi_{X}}
\arrow{ssw,l}{g} \node[2]{M(X)} \arrow{s,r}{\varphi_{X'}}
\arrow{ssw,r}{g}\\ \node[2]{X} \arrow[2]{e,t}{h} \arrow{ssw,l}{f}
\node[2]{X'} \arrow{ssw,l}{f'}\\ \node{M(Y)}
\arrow{s,l}{\varphi_{Y}} \arrow[2]{e,t}{\beta} \node[2]{M(Y)}
\arrow{s,l}{\varphi_{Y'}}\\ \node{Y} \arrow[2]{e,t}{k}
\node[2]{Y'}
\end{diagram}
\]
where $h,k$ are $\hat{t}$-equivalences, the vertical maps are all
$\mathbf{s}(f)$-maximal maps , the back and front face diagrams
are homotopy pushout diagrams and $$\alpha=\bigvee_i\alpha_i:M(X)
\to M(X),\beta=\bigvee_j\beta_j:M(Y) \to M(Y)$$ where
$$\alpha_i:M(H_i(X,\mathbb{Z})/torsion,i)\to
M(H_i(X,\mathbb{Z})/torsion,i)$$
$$\beta_j:M(H_j(Y,\mathbb{Z})/torsion,j)\to
M(H_j(Y,\mathbb{Z})/torsion,j)$$ and
$$M(X)=\sum_iM(H_i(X,\mathbb{Z})/torsion,i),M(Y)=\sum_jM(H_j(Y,\mathbb{Z})/torsion,j)$$
\end{thm}
\begin{proof}
Take $\hat{t}$-equivalences $h:X \to Y, k:Y \to Y'$  as granted by
Theorem \ref{T:dual2.2} such that $kf \sim f' h$. We will try to
define the other maps in the above diagram so that the Theorem is
true. As in the dual case , we can assume there is only one $l$
such that $H_i(X,\mathbb{Q})=H_i(Y,\mathbb{Q})=0$ if $i\neq l$.
Choose rational
co-H-space structures on $X' , Y'$ such that $h_{(0)},k_{(0)},f'$
are co-H-maps. As in the dual case we can choose elements
 $\{x_1,\cdots,x_{r_1}\} \in \pi_*(X)$,
 $\{x'_1,\cdots,x'_{r_2} \in \pi_*(Y)$ ,
 $\{y_1,\cdots,y_{r_1}\} \in \pi_*(X')$ and
 $\{y'_1,\cdots,y'_{r_2}\} \in \pi_*(Y')$ such that they project to
 bases of $\pi_*(X)\bigotimes \mathbb{Q}$ ,
 $\pi_*(Y)\bigotimes \mathbb{Q}$ ,
$\pi_*(X')\bigotimes \mathbb{Q}$ and $\pi_*(Y')\bigotimes \mathbb{Q}$
respectively. As in the dual case
 $\pi_*(f)/torsion$ is of the diagonal form with respect to these
bases and
\[
(h)_*(x_i)=\sum\lambda_{ii'} y_{i'} + v_i ,
(k)_*(x'_j)=\sum\lambda'_{jj'} y'_{j'} +w_j
\]
where $\det\lambda_{ii'},\det\lambda'_{jj'}$ are prime relative to
$\hat{t}$ and $v_i,w_j$ are torsions.
 Certainly we can find $\bar{v}_i$ and $\bar{w}_j$ so
that
\[
(h)_*(x_i)=\sum\lambda_{ii'} (y_{i'} + \bar{v}_{i'} ),
(k)_*(x'_j)=\sum\lambda'_{jj'} (y'_{j'} +\bar{w}_{j'})
\]
Let $\{\iota_1,\cdots,\iota_{r_1}\}$ be a base for
$H_*(M(X),\mathbb{Z})/torsion$ and
$\{\iota'_1,\cdots,\iota'_{r_2}\}$ be a basis for
$H_*(M(Y),\mathbb{Z})/torsion$. Define  $\varphi_Y,\varphi_X$  by
$s_l(f)x_i=(\varphi_X)_*(\iota_i)$,
$s_l(f)x'_j=(\varphi_Y)_*(\iota'_j)$ respectively and define $g$
so that its matrix with respect to
$\{\iota_1,\cdots,\iota_{r_1}\}$,
$\{\iota'_1,\cdots,\iota'_{r_2}\}$ is the same as that of
$\pi_*(f)/torsion$ . Let $\varphi_{X'}$ be the map which sends
each $\iota_i$ to $s_l(f)(y_i +\bar{v}_i)$ and $\varphi_{Y'}$ be
the map which sends each $\iota'_j$ to $s_l(f)(y'_j +\bar{w}_j)$
and set $\alpha,\beta$ to be the maps which send $\iota_i$ and
$\iota'_j$ to $\sum\lambda_{ii'} \iota_{i'}$ and
 $\sum\lambda'_{jj'} \iota'_{j'}$ , respectively. Obviously the only thing
we have to verify is that $\varphi_{Y'} g \sim f' \varphi_{X'} $.
It is easy to know from diagram chasing argument that $g
\varphi_{Y'}\alpha \sim f' \varphi_{X'} \alpha$. Now $(g
\varphi_{Y'} )_*(\iota_i)-(\varphi_{X'}f')_*(\iota_i)=u_i$ which
are torsions of order dividing $\hat{t}$ for all $i$. $u_i=0$ for
all $i$ iff $\sum \lambda_{ii'}u_{i'}=0$ for  all $i$ iff $g
\varphi_{Y'}\alpha \sim f' \varphi_{X'} \alpha$ since $\det
\lambda_{ii'}$ is prime relative to $\hat{t}$.
\end{proof}
\begin{proof}[Proof of Theorem\ref{T:main2}]
Without loss of generality, we can assume that there is an integer
$l$ such that $Q\pi_i(X)\bigotimes \mathbb{Q}=Q\pi_i(Y)\bigotimes \mathbb{Q}=0$ if
$i\neq l$. By Theorem \ref{T:dual3.4} and Theorem \ref{T:modifco}
it follows that there exists a surjection $\xi':T \to G(f)$ where
\[
T=\{(\alpha,\beta)|\beta g \sim g\alpha , \text{$\alpha$ and
$\beta$ are $\hat{t}$-equivalences}\} .
\]
In other words, there is a surjection $\xi':T' \to G(f)$ where
\[
T'=\{(A_1,A_2)|A_1 \in M(\mathbb{Z},r_i), (\det A_1,\hat{t} )=1,
i=1,2,  A_2C=CA_1\} .
\]
and $M(\mathbb{Z},n)$ is the set of $n\times n$ matrices and
$A_1,A_2,C$ represent $\alpha,\beta,g$ respectively. Given
$(A_1,A_2), (B_1,B_2)\in T^{'}$ such that $\det A_i \equiv \det
B_i $ (mod $\hat{t}$) for $i=1,2$, then an elementary calculation
shows that $A_1 =B_1 G_1+\hat{t}H_1$, $A_2= B_2 G_2 +\hat{t}H_2$ ,
$(G_1,G_2)\in T^{'} , \det G_1 \equiv \det G_2 \equiv 1(mod\hat{t})$ and
$G_2C=CH_2$. Claim \ref{T:claim} implies that we can assume $\det
G_1=\det G_2=1$. It follows from Theorem \ref{T:modifco} and
Theorem \ref{T:4.6-co} and Remark \ref{T:3.4re-co} that
$\xi^{'}(A_1,A_2)=\xi^{'}(B_1,B_2)$. Thus $\xi^{'}$ factors
through the map $\det:T^{'} \to
(\mathbb{Z}^*_{\hat{t}})^{l(X)+l(Y)}$. The rest of the proof is
dual to that of Hurvitz\cite{hurv}.
\end{proof}
\section{Some results about co-H-spaces}
In this section we will prove a result about co-H-spaces which is
needed in the last section and which may have independent
interest. This result concerns about the relative homology
decomposition of a map . Zabrodsky had proved that the Postnikov
decomposition of an H-map consists of H-spaces and H-maps (c.f.,
Corollary 2.3.2 in \cite{az1}).  The dual result is also true
under our assumption that all co-H-spaces we are concerned are
1-connected.
\begin{thm}\label{T:co-H}
Let $f:X \to Y$ be a co-H-map between co-H-spaces. Then the
homology decomposition of $f$ consists of co-H-spaces and
co-H-maps.
\end{thm}
\begin{rem}
When the map $f$ is the inclusion of the base point into the space
$X$ , the homology decomposition of $f$ is the Moore decomposition
of $X$. In that special case the above Theorem was already
obtained by M.Arkowitz\cite{arko2} and M.Golasinski and John R.
Klein\cite{klein}
\end{rem}
To prove the Theorem above we will follow Zabrodsky's approach to
the dual result. First let us recall results about the obstruction
of a map between co-H-spaces to be a co-H-map , see
Arkowitz\cite{arko1} for a comprehensive survey about co-H-space.

Let $XbX$ be the space of paths in $X\times X$ beginning in
$X\bigvee X$ and ending at the base point and let $i: XbX \to
X\bigvee X$ assign to a path its initial point. Then we have a
short exact sequence
\[
0 \to [\Omega X,\Omega(XbX)]\overset{(\Omega i)_*}{\to}[\Omega
X,\Omega(X\bigvee X)]\overset{(\Omega j)_*}{\to} [\Omega
X,\Omega(X\times X)] \to 0
\]
where $i_1,i_2$ are the inclusions to the first and second factors
of $X\bigvee X$. The comultiplication $\nu:X X \bigvee X$
determines an element $$\mu_X=-\Omega i_1 -\Omega i_2 +\Omega \nu
\in [\Omega X,\Omega(X\bigvee X)]$$ such that $(\Omega
j)_*(\mu_X)=0$. Thus there exists a unique element $H(\nu)\in
[\Omega X,\Omega(XbX)]$ such that $(\Omega i)_*(H(\nu))=\mu_X$. We
call $H(\nu)$ the dual Hopf construction (applied to $\nu$). It is
well known that $\beta:\Sigma A \to X $ is a co-H-map if and only
if $H(\nu)\bar{\beta}=0$  in $[A,\Omega(XbX)]$ where $\Sigma A$
has the standard comultiplication on the suspension and
$\bar{\beta}$ is the adjoint of the map $\beta$ .
\begin{lem}\label{T:obst}
Let $f:(X,\nu_X) \to (Y,\nu_Y)$ be a co-H-map between co-H-spaces
and $\beta:\Sigma A \to X$ be a map. Then
\[
H(\nu_Y)(\Omega f\bar{\beta})=\Omega(fbf)H(\nu_X)\bar{\beta}
\]
\end{lem}
\begin{proof}
It follows  from the equation $\mu_Y (\Omega
f\bar{\beta})=\Omega(fbf)\mu_X \bar{\beta}$ which can be verified
directly.
\end{proof}
It is well known that $XbX$ has the weak homotopy type of $\Sigma
\Omega X \bigwedge \Omega X$. The following is a special case of a
general result given by Golasinski and Klein.
\begin{lem}\cite{klein}\label{T:klein}
If $f:X \to Y$ is a co-H-map between co-H-spaces and $f$ is
$n$-connected with $n \geq 1$. Then $fbf$ is $n+1$-connected.
\end{lem}
\begin{rem}
The authors learned of this result through the discussion list of
Hopf archive. Thanks go to all those who have shared information
about this with us.
\end{rem}
\begin{prop}\label{T:prop}
Let $f:(X,\nu_X) \to (Y,\nu_Y)$ be a co-H-map. Suppose $f$ is
n-connected with $n\geq 2$. Given any map $\beta:\Sigma
M(G,n-1)=M(G,n) \to X$ so that $f\beta \sim *$, then $\beta$ is a
co-H-map.
\end{prop}
\begin{proof}
By Lemma \ref{T:obst} , we have
\[
*=H(\nu_Y)(*)=H(\nu_Y)(\Omega
f\bar{\beta})=\Omega(fbf)H(\nu_X)\bar{\beta}
\]
On the other hand , $H(\nu_X)\bar{\beta} \in
[M(G,n-1),\Omega(XbX)]=\pi_{n-1}(\Omega(XbX),G)$. Thus  $\beta$ is
a co-H-map if  $H(\nu_X)\bar{\beta}=0$. Now
$\Omega(fbf)H(\nu_X)\bar{\beta}=0$ . The universal coefficient
formula for homotopy group with coefficient and
Lemma \ref{T:klein} imply that $\Omega(fbf)$ induces a
monomorphism which gives the desired equation
$H(\nu_X)\bar{\beta}=0$.
\end{proof}
\begin{proof}[Proof of Theorem\ref{T:co-H}]
Now the proof of Theorem \ref{T:co-H} follows from  an induction by
the homology decomposition of the co-H-map $f$ together with
Proposition \ref{T:prop} and the fact that cofiber of a co-H-map is
 a co-H-space and the inclusion of the cofiber is a co-H-map.
\end{proof}

--------------------------------------------------------------

----------------------------------------------------------------

\end{document}